\documentclass{amsproc}

\newtheorem{theorem}{Theorem}[section]
\newtheorem{lemma}[theorem]{Lemma}

\theoremstyle{definition}

\newtheorem{example}[theorem]{Example}

\theoremstyle{remark}
\newtheorem{remark}[theorem]{Remark}

\newtheorem{corollary}[theorem]{Corollary}

\numberwithin{equation}{section}

\newcommand{\fm}{\mathfrak m}
\newcommand{\la}{\lambda}

\newcommand{\Ap}{\textrm{Ap}}
\newcommand{\w}{\omega}

\begin{document}

\title{Tangent cones of numerical semigroup rings }

\author{Teresa Cortadellas Ben\'{i}tez}
\address{Departament d'\`{A}lgebra i Geometria, Universitat de Barcelona, Gran Via 585, 08007 Barcelona}
 \email{terecortadellas@ub.edu}
\thanks{Partially supported by MTM2007-67493}

\author{Santiago Zarzuela Armengou}
\address{Departament d'\`{A}lgebra i Geometria, Universitat de Barcelona, Gran Via 585, 08007 Barcelona}
\email{szarzuela@ub.edu}
\thanks{}

\date{}



\begin{abstract}
In this paper we describe the structure of the
tangent cone of a numerical semigroup ring $A=k[[S]] \subseteq
k[[t]]$ with multiplicity $e$ (as a module over the Noether
normalization determined by the fiber cone of the ideal generated
by $t^e$) in terms of some classical invariants of the corresponding numerical semigroup. Explicit computations are also made by using the GAP system.
\end{abstract}

\maketitle

\section{Introduction}
Let $(A,\fm)$ be an one dimensional Cohen-Macaulay local ring with
infinite residue field, embedding dimension $b$, reduction number
$r$ and multiplicity $e$. Let $(x)$ be a minimal reduction of
$\fm$.

\medskip

In \cite{CZ1} and \cite{CZ2} the authors have observed that the Noether
normalization
\[ F_{\fm}(x):= \bigoplus_{n\geq 0} \frac{x^nA}{x^n\fm
}\hookrightarrow G(\fm):= \bigoplus_{n\geq 0}
\frac{\fm^n}{\fm^{n+1}}\]
 provides a decomposition of $G(\fm)$ as
direct sum of graded cyclic $F_{\fm}(x)$-modules of the form
\[G(\fm)\cong \bigoplus_{i=0}^{r}\left(F_{\fm}(x)(-i)\right)^{\alpha_i}
\bigoplus_{i=1}^{r-1}\bigoplus_{j=1}^{r-i-1}\left(\frac{F_{\fm}(x)}{(x^\ast)^jF_{\fm}(x)}(-i)\right)^{\alpha_{i,j}},
\] with $\alpha_0=1$, $\alpha_r\neq 0$ and
$\displaystyle{\sum_{i=1}^{r}\alpha_i}=e-1$.

\medskip

In the more general context of the study of fiber cone of ideals with analytic spread one, the authors analyze in \cite{CZ1} the information provided by the set of invariants $\{ \alpha_i,
\alpha_{i,j}\}$ (that we call the invariants of the tangent cone) in order to study, for instance, the Cohen-Macaulay or Buchsbaum properties of the tangent cone; whereas in \cite{CZ2} the connection between these invariants and the so called microinvariants introduced by Juan Elias \cite{E} in the geometric case, and other invariants introduced by Valentina Barucci and Ralf Fr\"{o}berg \cite{BF} in the numerical semigroup case is study in detail.

\medskip

The main purpose of this paper is to compute explicitly the values
of the set of invariants $\{ \alpha_i, \alpha_{i,j}\} $ for the
tangent cones of numerical semigroups rings in terms of the
invariants introduced by Barucci-Fr\"{o}berg in \cite {BF}. In
particular, these computations can be performed numerically by
using the GAP system \cite{GAP4}, as we show with several
examples. We note that for the case of the microinvariants Elias
himself has described the explicit computations for semigroup
rings in \cite[Section 4]{E}.

\medskip

The content of this work arose during the talk that, reporting the results in \cite{CZ2}, the second author gave at the Exploratory Workshop on Combinatorial Commutative Algebra and Computer Algebra held in Mangalia, Romania, in May 2008. Prof. J\"{u}rgen Herzog suggested during the talk that for the numerical semigroup case one should be able to compute explicitly the invariants of the tangent cone. The authors want to thank J. Herzog for this suggestion. Also, the second author would like to thank the organizers of the workshop, Alexandru Bobe, Viviana Ene, and Denis Ibadula from the Ovidius University in Constanta, for the invitation to participate in it and the excellent organization and warm atmosphere during the workshop.

\section{Tangent cones of numerical semigroups}

Let $\mathbb N$ be the set of non-negative integers.  Recall that a numerical
semigroup $S$ is a subset of $\mathbb N$ that is closed under
addition, contains the zero element and has finite complement in
$\mathbb N$.  The least integer belonging to $S$ is known as the
multiplicity of $S$ and it is denoted by $e(S)$. Reciprocally, the
greatest integer not belonging to $S$ is known as the Frobenius
number of $S$ and it is denoted by $F(S)$.

\medskip

A numerical semigroup $S$ is always finitely generated; that is, there
exist integers $n_1,\dots,n_l$ such that $S=\langle n_1,\dots
,n_l \rangle=\{ \alpha_1 n_1 +\cdots  +\alpha_l n_l ; \alpha_i\in
\mathbb N\}$. Moreover, every numerical semigroup has an unique minimal
system of generators $n_1,\dots ,n_{b(S)}$.

\medskip

A relative ideal of $S$ is a nonempty set $I$ of non-negative integers such
that $I+S\subset I$ and $d+I\subseteq S$ for some $d\in S$. An
ideal of $S$ is then a relative ideal of $S$ contained in $S$. If
$i_1,\dots ,i_k$ is a subset of non-negative integers, then the set
$\{i_1,\dots ,i_k\}+S=(i_1+S)\cup \cdots \cup (i_k+S)$ is a relative ideal of $S$ and
$i_1,\dots ,i_k$ is a system of generators of
$I$. Note that, if $I$ is an ideal of $S$, then $I\cup \{0\}$ is a
numerical semigroup and so $I$ is finitely generated. We denote by
$M$ the maximal ideal of $S$, that is, $M=S\setminus \{0\}$. $M$ is then
the ideal generated by a system of generators of $S$.
If $I$ and $J$ are relative ideals of $S$ then $I+J=\{i+j;i\in I, j\in
 J\}$ is also a relative ideal of $S$. Finally, we denote by $\Ap (I)$ the Apery
set of $I$ with respect to $e(S)$, defined as the set of the
smallest elements in $I$ in each residue class module $e(S)$.

\medskip

Let $V=k[[t]]$ be the formal power series ring over a field $k$.
Given a numerical semigroup $S=\langle n_1,\dots ,n_b \rangle$
minimally generated by $0<e=e(S)=n_1<\cdots < n_b=n_{b(S)}  $ we
consider the ring associated to $S$ defined as
$A=k[[S]]=k[[t^{n_1},\dots ,t^{n_b}]] \subseteq V$. Let
$\fm=(t^{n_1},\dots ,t^{n_b})$ be the maximal ideal ideal of $A$.
Then $A$ is a Cohen-Macaulay local ring of dimension one with
multiplicity $e$ and embedding dimension $b$. Moreover, the
conductor $(A:V)=t^{C}V$ with $C=F+1$ where $F=F(S)$ is the
Frobenius number of $S$. These kind of rings are known as numerical
semigroup rings. The ideals $(t^{i_1},\dots,t^{i_k})$ of $A$ are
such that for $v$, the $t$-adic valuation, $
v((t^{i_1},\dots,t^{i_k}))=\{i_1,\dots ,i_k\}+S$. In particular, for
the ideals $\fm^n$ one has $v(\fm^n)=nM = M + \stackrel{n}{\cdots} +M$.
Note that $(n+1)M\subseteq nM$ for $n\geq 0$ (we will set $\fm^0:=A$).

\medskip

Let $A=k[[S]]$ be a numerical semigroup ring of multiplicity $e$.
Then, the element $t^e$ generates a minimal reduction of $\fm$. In
terms of semigroups, $(n+1)M=e+nM$ for all $n\geq r$, the reduction number of $\fm$; that is,
$r$ is the smallest integer $n$ such that $\fm^{n+1}=t^e\fm^n$ (in our case the
reduction number does not depend on the minimal reduction).

\medskip

A crucial point for our results is the use of a fact proved by Barucci-Fr\"{o}berg \cite[Lemma 2.1]{BF} in the more general context of one-dimensional equicharacteristic analytically irreducible and residually rational doamins. For completenees, we give an easy proof of it for the particular case we deal with in this paper. Set
$k[[t^e]]=W\hookrightarrow A$.


\begin{lemma}
Let $I$ be an ideal of $S$ and $\mathfrak I$ the ideal of $A$ generated by $\{ t^n\}_{n\in
I}$. If $\Ap (I)=\{ \w_0,\dots ,\w_{e-1}\}$ is the Apery set of
$I$ with respect to $e$, then $\mathfrak I$ is a free $W$-module
generated by $t^{\w_0},\dots t^{\w_{e-1}}$.
\end{lemma}

\begin{proof}
Let $n\in I$. If $n\equiv i$ mod $e$ then $n=\w_i +\alpha e$ for
some $\alpha \geq 0$. So $t^n =(t^e)^{\alpha } t^{\w_i}\in
Wt^{\w_0}+\cdots +Wt^{\w_{e-1}}$ and $\mathfrak I
=Wt^{\w_0}+\cdots +Wt^{\w_{e-1}}$. Observe that the sum is direct
since in each summand the elements are monomials in $t$ with
exponents in different residue classes mod $e$.
\end{proof}

In particular, we may write the powers of the maximal ideal as a
direct sum of cyclic $W$-modules.

\begin{lemma}
\label{BF}For each $n\geq 0$ there exist non-negative integers $\omega_{n,0},
\dots , \omega_{n,e-1}$ such that
\[\fm^n
=Wt^{\omega_{n,0}}\oplus \cdots  \oplus Wt^{\omega_{n,e-1}},\]
 with
$\omega_{n+1,i}=\omega_{n,i}+e\cdot \epsilon$ and $\epsilon \in
\{0,1\} $.
\end{lemma}

\begin{proof} Observe first that if $\Ap(S)=\{ 0,\omega_1, \dots ,\omega_{e-1}\}$
is the  Apery set of $S$ (with respect to the multiplicity $e$),
then the Apery set of $M=S\setminus \{0\}$ the  maximal ideal of
$S$ is $\Ap(M)=\{e,\omega_1, \dots ,\omega_{e-1}\}$.

Now, for each $n\geq 1$ let $\Ap (nM)=\{ \omega_{n,0},\dots
,\omega_{n,e-1}\}$ be the Apery set of $nM$. If  $\omega_{n,i}\in
(n+1)M$ then $\omega_{n+1,i}=\omega_{n,i}$. Otherwise
$\omega_{n,i}+e\in (n+1)M$ and belongs to the same residue
class of $\omega_{n,i}$ module $e$. Since $\omega_{n+1,i}\leq
\omega_{n,i}+e$ by definition, it follows that $\omega_{n,i}+e-
\omega_{n+1,i}=\alpha e$ for some $\alpha \geq 1$. On the other
hand, if $\alpha \geq 2$ then $\omega_{n,i}=
\omega_{n+1,i}+(\alpha -1)e \in (n+1)M$ which contradicts the
assumption, so $\alpha = 1$. Now the proof is concluded by
applying the above lemma to $nM$ for $n\geq 0$ (where $0M:=S$).
\end{proof}

Observe that for each $n\geq 0$ and each $0\leq i \leq e-1$, $Wt^{\omega_{n+1,i}} \subseteq Wt^{\omega_{n,i}}$. Also, that for $n\geq r$ we have
$\omega_{n+1,i}=\omega_{n,i}+e$.

\medskip

Our next result gives a description of the set of invariants $\{ \alpha_i,
\alpha_{i,j}\}$ of the tangent cone in terms of the Apery sets of the family of ideals $nM$, for $0\leq n \leq r$. Previously, and just for the purposes of this paper, we introduce the following notation:

\medskip

Let $E = \{a_0, \dots , a_n\}$ be a set of integers. We call it a
ladder if $a_0 \leq \cdots \leq a_n$. Given a ladder, we say that
a subset $L = \{ a_i, \dots, a_{i+k}\}$ with $k\geq 1$ is a
landing of length $k$ if $a_{i-1} < a_i = \cdots = a_{i+k} <
a_{i+k+1}$ (where $a_{-1} = -\infty$ and $a_{n+1} = \infty$). In
this case, the index $i$ is the beginning of the landing: $s(L)$
and the index $i+k$ is the end of the landing: $e(L)$. A landing
$L$ is said to be a true landing if $s(L) \geq 1$. Given two
landings $L$ and $L'$, we set $L < L'$ if $s(L) < s(L')$. Let
$l(E)+1$ be the number of landings and assume that $L_0< \cdots
<L_{l(E)}$ is the set of landings. Then, we define following
numbers:

\begin{itemize}

\item[$\cdot$] $s_j(E)= s(L_j)$, $e_j(E) = e(L_j)$, for each $0\leq j \leq l(E)$;

\medskip

\item[$\cdot$] $c_j(E) = s_j - e_{j-1}$, for each $1 \leq j \leq l(E)$.

\end{itemize}

\noindent Note that the above numbers are defined under the conditions

\medskip

\[
\begin{split}
& a_0 = \cdots = a_{e_0(E)} < \cdots < a_{e_0(E) + c_1(E)} = \cdots = a_{e_1(E)} < \\
& \cdots  \cdots  \\
& < a_{e_{(l(E)-1)} + c_{l(E)}} = \cdots = a_{e_{l(E)}} < a_{e_{l(E)}+1} < \cdots < a_{n}
\end{split}
\]

\medskip

\begin{theorem}
\label{main} Let $A=k[[S]]$ be a numerical semigroup ring of
multiplicity $e$ and reduction number $r$. Let $M$ be the maximal
ideal of $S$ and put \[ \Ap (nM)=\{\w_{n,0},\dots ,\w_{n,i}, \dots
,\w_{n,e-1}\}\] for $0\leq n \leq r$.

For any $1\leq i \leq e-1$, consider the ladder of values
$W^i=\{\w_{n,i}\}_{0\leq n\leq r}$ and define the following integers:
\begin{enumerate}
\item $l_i = l(W^i)$;
\item $d_i = e_{l_i}(W^i)$;
\item $b_j^i=e_{j-1}(W^i)$ and $c_j^i=c_j(W^i)$, for $1\leq j\leq l_i$.
\end{enumerate}
Then
\[G(\fm)\cong F_{\fm}(t^e)\oplus \bigoplus_{i=1}^{e-1}\left(
F_{\fm}(t^e)(-d_i)
 \bigoplus_{j=1}^{l_i} \frac{F_{\fm}(t^e)}{((t^e)^{\ast})^{c_j^i}F_{\fm}(t^e)}
(-b_j^i)\right).\]

\end{theorem}

\begin{proof}
For all $n\geq 0$, we have by Lemma \ref{BF} that $\fm^n=
Wt^{\omega_{n,0}}\oplus\dots \oplus Wt^{\omega_{n,e-1}}$, with
$\w_{0,n}=t^{en}$ for all $n$, and so we have the following
commutative diagram of graded rings
\[
\begin{array}{ccccc}
G(\fm) &\cong &\displaystyle{\bigoplus_{n\geq 0}\left(
\bigoplus_{i=0}^{e-1}\frac{Wt^{\w_{n,i}}}{Wt^{\w_{n+1,i}}}\right)} & = &G\\
\uparrow & & \uparrow \\ F_{\fm}(t^e) &\cong
&\displaystyle{\bigoplus_{n\geq 0}\frac{(t^e)^n W }{(t^e)^{n+1}W}}
& = &F
\end{array}
\] and we can read the structure of $G(\fm)$ as $F_´{\fm}(t^e)$-module
as the structure of $G$ as $F$-module. Note that $G$ may also be written as
$$G = \bigoplus_{i=0}^{e-1}\left(
\bigoplus_{n\geq 0}\frac{Wt^{\w_{n,i}}}{Wt^{\w_{n+1,i}}}\right)$$

Now, let us fix $1\leq i \leq e-1$. Assume first that $b_1^i=d_i$. Then, the component of degree
degree $n$ of $\displaystyle{\bigoplus_{n\geq 0}
\frac{Wt^{\w_{n,i}}}{Wt^{\w_{n+1,i}}}}$ is

\[ \left[ \bigoplus_{n\geq 0}\frac{Wt^{\w_{n,i}}}{Wt^{\w_{n+1,i}}}
\right]_n =
\begin{cases} \begin{array}{ll} 0 &\textrm{ if } 0\leq n <d_i \\
\displaystyle{\frac{Wt^{\w_{0,i}+e(n-d_i)}}{Wt^{\w_{0,i}+e(n-d_i+1)}}}&
\textrm{ if } n\geq d_i\end{array}\end{cases}$$ with
$$\frac{Wt^{\w_{0,i}+e(n-d_i)}}{Wt^{\w_{0,i}+e(n-d_i+1)}}=
\frac{Wt^{\w_{0,i}}\cdot t^{e(n-d_i)}}{Wt^{\w_{0,i}}\cdot
t^{e(n-d_i+1)}} \cong \frac{W(t^e)^{(n-d_i)}}{W(t^e)^{(n-d_i+1)}}
\]
 and
\[\bigoplus_{n\geq 0} \frac{Wt^{\w_{n,i}}}{Wt^{\w_{n+1,i}}}\cong
F(-d_i).\]

Assume now that  $b_1^i<d_i$. Then the piece of degree $n$ of
$\displaystyle{\bigoplus_{n\geq 0} \frac{Wt^{\w_{n,i}}}{Wt^{\w_{n+1,i}}}}$
is in this case

\[ \left[ \bigoplus_{n\geq 0}\frac{Wt^{\w_{n,i}}}{Wt^{\w_{n+1,i}}}
\right]_n =
\begin{cases} \begin{array}{ll} 0 &\textrm{ if } 0\leq n <b^i_1 \\
\displaystyle{\frac{Wt^{\w_{0,i}+e(n-b^i_1)}}{Wt^{\w_{0,i}+e(n-b^i_1+1)}}}&
\textrm{ if } b^i_1\leq n<  b^i_1+c^i_1 \\ &\vdots \\
 0 &\textrm{ if } b^i_{j-1}+ c^i_{j-1}\leq n <b^i_j \\
\displaystyle{\frac{Wt^{\w_{0,i}+e(n+c^i_1+\cdots
+c^i_{j-1}-b^i_j)}}{Wt^{\w_{0,i}+e(n+1+c^i_1+\cdots
+c^i_{j-1}-b^i_j)}}}& \textrm{ if } b^i_j\leq n< b^i_j+c^i_j \\
&\vdots  \\
 0 &\textrm{ if } b^i_{l_i}+ c^i_{l_i}\leq n <d_i \\
\displaystyle{\frac{Wt^{\w_{0,i}+e(n+c^i_1+\cdots
+c^i_{l_i}-d_i)}}{Wt^{w_{0,i}+e(n+1+c^i_1+\cdots
+c^i_{l_i}-d_i)}}}& \textrm{ if }n\geq d_i \\
\end{array}\end{cases}\]

\noindent
with $W$-isomorphisms
\[\frac{Wt^{\w_{0,i}+e(n+c^i_1+\cdots
+c^i_{j-1}-b^i_j)}}{Wt^{\w_{0,i}+e(n+c^i_1+\cdots
+c^i_{j-1}-b^i_j+1)}} \cong
\frac{W(t^e)^{(n-b^i_j)}}{W(t^e)^{(n-b^i_j+1)}} \]
 for $1\leq j \leq l_i$ and $b^i_j\leq n\leq b^i_j+c^i_j$,
\[
 \frac{Wt^{\w_{0,i}+e(n+c^i_1+\cdots
+c^i_{l_i}-d_i)}}{Wt^{\w_{0,i}+e(n+1+c^i_1+\cdots
+c^i_{l_i}-d_i)}}\cong \frac{W(t^e)^{(n-d_i)}}{W(t^e)^{(n-d_i)}}\]
for $n\geq d_i$ and
\[
\begin{split}
 \bigoplus_{n\geq 0}
\frac{Wt^{\w_{n,i}}}{Wt^{\w_{n+1,i}}}
 & = \bigoplus_{j=0}^{l_i} F\cdot (t^{\w_{0,i}+e(c^i_1+\cdots
+c^i_{j})})^{\ast} \\
&\cong  \bigoplus_{j=1}^{l_i}
\frac{F}{((t^e)^{\ast})^{c_j^i}F}(-b_j^i) \oplus F(-d_i).
\end{split}
\]
\end{proof}

\begin{example}
\label{5,6,13}
 Let $S=\langle 5,6,13 \rangle$.  It is easy to
prove that the maximal ideal of $S$ has reduction number 4 and
also to calculate the Apery sets of the ideals $nM$ for $n\geq 0$.
The following table shows these values for $n\leq 4$:

\begin{table}[ht]
\renewcommand\arraystretch{1.5}
\[
\begin{array}{|l|c|c|c|c|c|}\hline
\Ap(S)& 0&6&12&13&19\\ \hline \Ap(M)&5&6&12&13&19\\ \hline
\Ap(2M)&10&11&12&18&19\\ \hline
 \Ap(3M)&15&16&17&18&24\\ \hline
 \Ap(4M)&20&21&22&23&24 \\ \hline
 \end{array}
\]
\end{table}

\medskip

Then the tangent cone  $G$ of $k[[t^5,t^6,t^{13}]] $ has the
following structure over $F$ the fiber cone of $(t^5)$
\[
  F\oplus
F\cdot(t^6)^{\ast}\oplus F\cdot (t^{12})^{\ast}\oplus F\cdot
(t^{13})^{\ast}\oplus F\cdot (t^{18})^{\ast}\oplus F\cdot
(t^{19})^{\ast}\oplus F\cdot (t^{24})^{\ast},
\]
 and so isomorphic to \[
 F\oplus F(-1)\oplus F(-2) \oplus (F/(t^5)^{\ast}F)(-1)\oplus
F(-3)\oplus (F/(t^5)^{\ast}F)(-2)\oplus F(-4).\]
\end{example}

\begin{remark}
\label{table} Observe that the necessary information to determine
the structure of $G$ as $F$-module is contained in the table
\begin{table}[ht]
\renewcommand\arraystretch{1.5}
\[
\begin{array}{|c|c|c|c|c|c|c|}\hline
 \Ap(S)& \w_{0,0} &\w_{0,1}&\cdots  &
 \w_{0,i} &\cdots & \w_{0,e-1} \\
\hline
\Ap(M)& \w_{1,0} &\w_{1,1}&\cdots & \w_{1,i} & \cdots & \w_{1,e-1}\\
\hline \vdots& \vdots & \vdots &\vdots & \vdots & \vdots & \vdots
\\ \hline
\Ap(nM)&\w_{n,0} &\w_{n,1}&\cdots & \w_{n,i} & \cdots &
\w_{n,e-1}\\
\hline
\vdots & \vdots & \vdots &\vdots & \vdots & \vdots & \vdots \\
\hline
\Ap(rM) & \w_{r,0} &\w_{r,1}& \cdots & \w_{r,i} & \cdots  & \w_{r,e-1}\\
\hline  \end{array} \]
 \end{table}

\noindent that we call the Apery table of $S$.

 \medskip

 Thus, if we analyze the increment of the values by columns then we obtain  the
 values of $d_i$, $b^i_{j},\dots  ,c^i_{j}$ for $1\leq i\leq e-1$
 and $1\leq j\leq l_i$ and we may write, putting $x=t^e$,
\[G\cong F\oplus \bigoplus_{i=1}^{e-1}\left( F(-d_i)
 \bigoplus_{j=1}^{l_i} \frac{F}{(x^{\ast})^{c_j^i}F}
(-b_j^i)\right).\]
 Also, if we separate free and torsion submodules and collect the summands
  by the degrees of the generators we can rewrite the above
  expression in the form
\[ G\cong \bigoplus_{i=0}^{r}\left( F(-i)\right)^{\alpha_i}
\bigoplus_{i=1}^{r-1}\bigoplus_{j=1}^{r-i-1}\left(\frac{F}{(x^\ast)^jF}(-i)\right)^{\alpha_{i,j}}.
\]

Hence setting
$\beta_{0,i}=\alpha_i+\displaystyle{\sum_{j=1}^{r-i-1}\alpha_{i,j}}$
and $\beta_{1,i}=\displaystyle{\sum_{k+l=i}\alpha_{k,l}} $ we get
that
$$ 0\longrightarrow \bigoplus_{i=1}^{r-1}F(-i)^{\beta_{1,i}} \longrightarrow \bigoplus_{i=0}^{r-1} F(-i)^{\beta_{0,i}}
\longrightarrow 0
$$
gives a minimal graded free resolution, or equivalently the graded
Betti numbers of $G$ as $F$-module.
\end{remark}

\section{Computing the invariants of the tangent cone: Examples}

The GAP - Groups, Algorithms, Programming - is a system for
Computational Discrete Algebra \cite{GAP4}. On the basis of GAP, Manuel Delgado,
Pedro A. Garcia-S\'{a}nchez and  Jos\'{e} Morais have developed the
NumericalSgps package \cite{NumericalSgps}. Its aim is to
make available a computational tool to deal with numerical
semigroups. By Theorem \ref{main},
we can determine the structure of the tangent cone of a numerical
ring $k[[S]]\subset k[[t]]$ of multiplicity $e$ as a module over
the fiber cone of $(t^e)$ if we know the Apery sets of the sum
ideals $nM$, where $M$ is the maximal ideal of $S$. On the other hand, from the definition we have that the Apery set of $nM$ can
be calculated as
$$ \Ap (nM)=nM \setminus ((e+S)+ nM)$$

\noindent (see also \cite[Lemma 2.1 (2)]{BF}), a computation that can be performed by using the NumericalSgps package. The following examples are just a sample of these computations.

\begin{example}
\label{10,11,19}
 Let $S=\langle 10,11,19 \rangle$.  By using the
NumericalSgps package of GAP we calculate the  reduction number of
$M$ which is $8$ and also  the Apery sets of the ideals $nM$ for
$n\leq 8$. The following is the Apery table in this case:

\begin{table}[ht]
\renewcommand\arraystretch{1.5}
\[
\begin{array}{|l|c|c|c|c|c|c|c|c|c|c|}\hline
Ap(S)& 0&11&22&33&44&55&66&57&38&19\\ \hline
Ap(M)&10&11&22&33&44&55&66&57&38&19\\ \hline
Ap(2M)&15&21&22&33&44&55&66&57&38&29\\ \hline
Ap(3M)&20&31&32&33&44&55&66&57&48&39\\ \hline
Ap(4M)&25&41&42&43&44&55&66&67&58&49\\ \hline
Ap(5M)&30&51&52&53&54&55&66&77&68&59\\ \hline
Ap(6M)&35&61&62&63&64&65&66&77&78&69\\ \hline
Ap(7M)&40&71&72&73&74&75&76&77&88&79\\ \hline
Ap(8M)&45&81&82&83&84&85&86&87&88&89\\ \hline
\end{array}
\]
\end{table}
By Theorem \ref{main}, the tangent cone  $G$ of
$k[[t^{10},t^{11},t^{19}]] $ has the following structure over $F$
the fiber cone of $(t^{10})$:
$$F\oplus (F(-1))^2 \oplus F(-2)\oplus F(-3)\oplus F(-4)\oplus
F(-5)\oplus F(-6)\oplus F(-7)\oplus F(-8)
 \oplus$$$$(F/x^2F)(-3)\oplus (F/x^5F)(-2)$$
where $x:=(t^{10})^{\ast}$. Thus we get that the minimal graded free
resolution of $G$ as
 $F$-module has the following values for its Betti numbers:
\[
\begin{split}
&\beta_{0,1}=\beta_{0,2}=\beta_{0,3}=2,\\
&\beta_{0,0}=\beta_{0,4}=\beta_{0,5}=\beta_{0,6}=\beta_{0,7}=\beta_{0,8}=1,
\\
&\beta_{1,5}=\beta_{1,7}=1.
\end{split}
 \]
\end{example}

\begin{example}
\label{10,19,47}
 Let $S=\langle 10,19,47 \rangle$.  By using the
NumericalSgps package of GAP we calculate the  reduction number of
$M$ which is $9$ and also the Apery sets of the ideals $nM$ for
$n\leq 9$ and we get the following Apery table:

\begin{table}[ht]
\renewcommand\arraystretch{1.5}
\[
\begin{array}{|l|c|c|c|c|c|c|c|c|c|c|c|}\hline
Ap(S)& 0&141&132&113&94&85&66&47&38&19\\ \hline
Ap(M)&10&141&132&113&94&85&66&47&38&19\\ \hline
Ap(2M)&20&141&132&113&94&85&66&57&38&29\\ \hline
Ap(3M)&30&141&132&113&104&85&76&57&48&39\\ \hline
Ap(4M)&40&151&132&123&104&95&76&67&58&49\\ \hline
Ap(5M)&50&151&142&123&114&95&86&77&68&59\\ \hline
Ap(6M)&60&161&142&133&114&105&96&87&78&69\\ \hline
Ap(7M)&70&161&152&133&124&115&106&97&88&79\\ \hline
Ap(8M)&80&171&152&143&134&125&116&107&98&89\\ \hline
Ap(9M)&90&171&162&153&144&135&126&117&108&99\\ \hline
\end{array}
\]
\end{table}

As a consequence, the tangent cone  $G$ of $k[[t^{10},t^{19},t^{47}]] $ has
the following invariants as a module over  $F$ the fiber cone of
$(t^{10})$:
\[
\begin{split}
& \alpha_i=1 \textrm{ for } 0\leq i \leq 9, \\
& \alpha_{1,1}=\alpha_{6,1}= \alpha_{7,1}=1,
 \alpha_{2,1}=\alpha_{4,1}=\alpha_{5,1}=2,
  \alpha_{3,1}=3.
 \end{split} \]
And the graded Betti numbers of $G$ as  $F$-module are
\[
\begin{split}
 &\beta_{0,0}=1, \beta_{0,1}=2, \beta_{0,2}=3,\beta_{0,3}=4, \beta_{0,4}=3,
 \beta_{0,5}=3, \beta_{0,6}=2, \beta_{0,7}=2, \\
 &\beta_{1,2}=1,
\beta_{1,3}=2, \beta_{1,4}=3, \beta_{1,5}=2, \beta_{1,6}=2,
\beta_{1,7}=1, \beta_{1,8}=1.
\end{split}
\]
\end{example}

\section{Buchsbaum property of numerical semigroup rings}
In this section we analyze the Buchsbaum property of the tangent
cones of numerical semigroup rings.

\medskip

Let $A=k[[S]]\subseteq k[[t]]$ be  a numerical semigroup ring of
multiplicity $e$, embedding dimension $b$ and reduction number
$r$. Set $x=t^e$. Recall that $e= \mu(\fm^n)+\la(\fm^{n+1}/x\fm^n)$ for all $n\geq 0$ and so $b = \mu(\fm) \leq e$. Also, that $\mu(\fm^n)\geq n+1$ for $0\leq n \leq r$ and so $r \leq \mu(\fm^r) - 1 = e -1$. Let $G$ be
the tangent cone of $A$ and $F$ the fiber cone associated to the
ideal $(x)$. Assume that
\[ G\cong F\oplus
\bigoplus_{i=1}^{e-1}\left( F(-d_i)
 \bigoplus_{j=1}^{l_i} \frac{F}{(x^{\ast})^{c_j^i}F}
(-b_j^i)\right)\]
 where the integers $\{d_i,l_i, c_j^i, b_j^i;\,
1\leq i\leq e-1,\, 0\leq j \leq l_i\}$ are as in Theorem
\ref{main}.

\medskip

The tangent cone $G$ is Cohen-Macaulay if and only if $G$ is a
free graded module over $F$. In this case the structure of $G$ as
$F$-graded  module is
\[
 G\cong F\oplus
\bigoplus_{i=1}^{e-1}\left( F(-d_i)\right).
\]

\noindent Observe that this is equivalent to the fact that there are no true landings in the ladders determined by the columns of the the Apery table of $S$.

\begin{example}
\label{10,17,22,28} Let $S=\langle 10,17,22,28 \rangle $. For this
numerical semigroup we use GAP to calculate the reduction number (which is $4$)
and the Apery table of $S$:

\begin{table}[ht]
\renewcommand\arraystretch{1.5}
\[
\begin{array}{|l|c|c|c|c|c|c|c|c|c|c|}\hline
\Ap(S)& 0&51&22&73&34&45&56&17&28&39\\
\hline \Ap(M) & 10&51&22&73&34&45&56&17&28&39\\
 \hline \Ap(2M)&  20&51&32&73&34&45&56&27&38&39\\ \hline
 \Ap(3M)& 30&51&42&73&44&55&56&37&48&49\\ \hline
  \Ap(4M) & 40&61&52&73&54&65&66&47&58&59
  \\ \hline
 \end{array}
 \]
\end{table}

\noindent and so
\[
 G\cong F\oplus F(-1)^3 \oplus F(-2)^3  \oplus F(-3)^2
\oplus F(-4).
\]
\end{example}

It is well known that the tangent cone of a ring with reduction
number at most 2 is Cohen-Macaulay. This is obvious in our case from the fact that there is no room for possible landings in the associated Apery table. The following examples show
the structure of these tangent cones in the case of numerical
semigroup rings in terms of the values of the associated numerical
semigroup by using the main theorem of section 2.

\begin{corollary}
\label{r=1}
 Let $S$ be a numerical semigroup with reduction number $1$ and set $S=\langle e=n_0, \dots,
n_{b(S)-1} \rangle $. Then $b(S) = e$ and
\[
\begin{split}
G&=F\oplus F\cdot (t^{n_1})^{\ast} \oplus  \cdots
\oplus F\cdot (t^{n_{e-1}})^{\ast} \\
& \cong F\oplus  F(-1)^{e-1}
\end{split}
\]
\end{corollary}

\begin{proof}
Observe first that $r = 1$ if and only if $b(S) = \mu(\fm) = e$ (that is, $A$ is of minimal multiplicity). Hence the values $n_i$, $n_j$ must belong to different residue classes module $e$ for
$i\neq j$ and we may assume that $n_i \equiv i$ module $e$. Then,
$\Ap(S)=\{ 0, n_1, \dots ,n_{e-1}\}$ and $\Ap(M)=\{ e, n_1, \dots
,n_{e-1}\}$, which implies that $d_i=b_1^i=1$ for $1\leq i \leq e-1$.
\end{proof}

\begin{corollary}
\label{r=2}
 Let $S$ be a numerical semigroup with
reduction number $2$. Then
\[
 G\cong F\oplus  F(-1)^{b-1}\oplus F(-2)^{e-b}
\]
\end{corollary}
\begin{proof}
In this case there exist $\w_{i_j}$ for $1\leq j \leq b-1$ such
$S=\langle e,\w_{i_1},\dots, \w_{i_{b-1}} \rangle $ and  $M=\{e,\w_{i_1},\dots
,\w_{i_{b-1}}\}+S$. Moreover, the $\w_{i_j}$'s are not in $2M$ and they are in different residue classes mod. $e$.
Thus, the corresponding Apery table is

\begin{table}[ht]
\renewcommand\arraystretch{1.5}
\[
\begin{array}{|l|c|c|c|c|c|c|c|c|}\hline
\Ap(S)& 0&\w_1&\cdots &\w_{i_1}&\cdots&\w_{i_{b-1}}&\cdots &\w_{e-1} \\
\hline \Ap(M)& e&\w_1&\cdots &\w_{i_1}&\cdots&\w_{i_{b-1}}&\cdots
&\w_{e-1} \\ \hline \Ap(2M)& 2e&\w_1&\cdots
&\w_{i_1}+e&\cdots&\w_{i_{b-1}}+e&\cdots &\w_{e-1} \\ \hline
 \end{array}
 \]
 \end{table}

\noindent  which gives the structure of $G$ by theorem \ref{main}.
\end{proof}

If the multiplicity of $A$ is less or equal to $3$ then the reduction number is at most $2$. So the following is one of the next cases:

\begin{corollary} Let $S$ be a numerical semigroup of multiplicity $4$ and  embedding dimension $b$.
\begin{enumerate}
 \item If $b=4$ then $G \cong F\oplus  F(-1)^{3}$.
\item If $b=3$ then $r=2$ or $r=3$ and
\begin{enumerate}
\item $G \cong F\oplus  F(-1)^{2}\oplus F(-2)$ if $r=2$,
 \item $G \cong F\oplus  F(-1)\oplus F(-2)\oplus F(-3)\oplus
(F/(t^4)^{\ast}F)(-1)$ if $r=3$.
\end{enumerate}
 \item If $b=2$ then $G \cong F\oplus  F(-1)\oplus F(-2)\oplus
F(-3)$.
\end{enumerate}
\end{corollary}

\begin{proof}
We have that $1\leq r\leq 3$ and $2\leq b\leq 4$. Hence, it suffices to determine all the possible Apery tables in each case an then apply theorem \ref{main}.

Assume first that $b=4$. Then, $r=1$ and the result follows from lemma \ref{r=1}. Moreover, if $S=(4,\w_1,\w_2,\w_3)$ the Apery table is in this case

\begin{table}[ht]
\renewcommand\arraystretch{1.5}
\[
\begin{array}{|l|c|c|c|c|}\hline \Ap(S)& 0&\w_1&\w_2&\w_3\\ \hline
 \Ap(M)& 4&\w_1&\w_2&\w_3\\ \hline
 \end{array}
 \]
 \end{table}

Assume now that $b=3$ and  set $S=(4,\w_1,\w_2)$ Then, $\la(\fm^2/x\fm)=1$
and so $r\geq 2$. If $r=2$, equivalently, if $\mu(\fm^2)=4$, there
exists $\w_3 \in S$ such that the Apery table (after a possible permutation of the columns) is

\begin{table}[ht]
\renewcommand\arraystretch{1.5}
\[
\begin{array}{|l|c|c|c|c|}\hline \Ap(S)& 0&\w_1&\w_2&\w_3\\ \hline
 \Ap(M)& 4&\w_1&\w_2&\w_3\\ \hline
 \Ap(2M)& 8&\w_1+4&\w_2+4&\w_3\\ \hline
 \end{array}
 \]
 \end{table}

\noindent Otherwise $r=3$, equivalently, $\mu(\fm^2)=3$. Taking
lengths in the exact sequence
\[
 0\longrightarrow (\fm^3+x\fm)/x\fm \longrightarrow \fm^2/x\fm
\longrightarrow \fm^2/(\fm^3+x\fm)\longrightarrow 0
\]
we get that $\fm^3 \subseteq x\fm$. Hence, there exist $\w_1,
\w_2, \w_3 \in S$ such that the Apery table (after a possible permutation of the columns) is

\begin{table}[ht]
\renewcommand\arraystretch{1.5}
\[
\begin{array}{|l|c|c|c|c|}\hline \Ap(S)& 0&\w_1&\w_2&\w_3\\ \hline
 \Ap(M)& 4&\w_1&\w_2&\w_3\\ \hline
 \Ap(2M)& 8&\w_1+4&\w_2+4&\w_3\\ \hline
\Ap(3M)& 12&\w_1+4&\w_2+8&\w_3+4\\ \hline
 \end{array}
 \]
 \end{table}

\pagebreak

Finally, assume that $b=2$. Then $r=3$, $\mu (\fm^2)=3$, $\mu(\fm^3)=4$ and $\la(\fm^3/x\fm^2)=1$. Thus, taking
lengths in the exact sequence
\[0\longrightarrow \fm^3/x\fm^2 \longrightarrow x\fm/x\fm^2
\longrightarrow x\fm/\fm^3\longrightarrow 0
\]
we get that $\la(x\fm/\fm^3)=1$  and so $\fm^3$ is not
contained in $x\fm$. Now, the Apery table (after a possible permutation of the columns) is given by

\begin{table}[ht]
\renewcommand\arraystretch{1.5}
\[
\begin{array}{|l|c|c|c|c|}\hline \Ap(S)& 0&\w_1&\w_2&\w_3\\ \hline
 \Ap(M)& 4&\w_1&\w_2&\w_3\\ \hline
 \Ap(2M)& 8&\w_1+4&\w_2&\w_3\\ \hline
\Ap(3M)& 12&\w_1+8&\w_2+4&\w_3\\ \hline
 \end{array}
 \]
 \end{table}

\end{proof}

\begin{example}
\label{e=4} This example illustrates the above corollary. In each
case we give the specific Apery table associated to the semigroup.

\begin{enumerate}
\item Let $S=\langle 4,10,11,17 \rangle$. Then,

\begin{table}[ht]
\renewcommand\arraystretch{1.5}
\[
\begin{array}{|l|c|c|c|c|}\hline \Ap(S)& 0&17&10&11\\ \hline
 \Ap(M)& 4&17&10&11\\ \hline
 \end{array}
 \]
\end{table}

\noindent and
\[
G= F\oplus F\cdot (t^{10})^{\ast} \oplus F\cdot
(t^{11})^{\ast}F\cdot (t^{17})^{\ast}\cong F\oplus F(-1)^{3}.\]

\item Let $S=\langle 4,10,11 \rangle $. Then,
\begin{table}[ht]
\renewcommand\arraystretch{1.5}
\[
\begin{array}{|l|c|c|c|c|}\hline \Ap(S)& 0&21&10&11\\ \hline
 \Ap(M)& 4&21&10&11\\ \hline
 \Ap(2M)& 8&21&14&15\\ \hline
 \end{array}
 \]
 \end{table}

\noindent and
\[
G= F\oplus F\cdot (t^{10})^{\ast} \oplus F\cdot
(t^{11})^{\ast}F\cdot (t^{21})^{\ast}\cong F\oplus F(-1)^{2}\oplus
F(-2).\]

\item Let $S=\langle 4,11,29 \rangle $. Then,
\begin{table}[ht]
\renewcommand\arraystretch{1.5}
\[
\begin{array}{|l|c|c|c|c|}\hline \Ap(S)& 0&29&22&11\\ \hline
 \Ap(M)& 4&29&22&11\\ \hline
 \Ap(2M)& 8&33&22&15\\ \hline
\Ap(3M)& 12&33&26&19\\ \hline
 \end{array}
 \]
 \end{table}

\noindent and
\[
\begin{split}
G &= F\oplus F\cdot (t^{11})^{\ast} \oplus F\cdot
(t^{22})^{\ast}F\cdot (t^{33})^{\ast} \oplus F\cdot
(t^{29})^{\ast}\\ &\cong F\oplus F(-1)\oplus F(-2)\oplus F(-3)
\oplus (F/(t^4)^{\ast}F)(-1). \end{split}
 \]

\item Let $S=\langle 4,11 \rangle $. Then,
\begin{table}[ht]
\renewcommand\arraystretch{1.5}
\[
\begin{array}{|l|c|c|c|c|}\hline \Ap(S)& 0&33&22&11\\ \hline
 \Ap(M)& 4&33&22&11\\ \hline
 \Ap(2M)& 8&33&22&15\\ \hline
\Ap(3M)& 12&33&26&19\\ \hline
 \end{array}
 \]
 \end{table}

\noindent and
\[
 G = F\oplus F\cdot (t^{11})^{\ast} \oplus F\cdot
(t^{22})^{\ast}F\cdot (t^{33})^{\ast}\cong F\oplus F(-1)\oplus
F(-2)\oplus F(-3).
 \]

\end{enumerate}
\end{example}

The tangent cone $G$ is Buchsbaum if and only if $G_{+}\cdot
H^0_{G_{+}}(G) =0$. Moreover, as observed in \cite{CZ1}, $H^0_{G_{+}}(G)$ coincides with $T(G)$,
the $F$-torsion submodule of $G$. As a consequence, if $G$ is Buchsbaum there cannot exist elements of order $>1$ in $T(G)$ and
then $c^i_j=1$ for al $i$ and $j$. That is, if $G$ is Buchsbaum
then
\[
 G\cong F\oplus \bigoplus_{i=1}^{e-1}\left( F(-d_i)
 \bigoplus_{j=1}^{l_i} \frac{F}{(x^{\ast})F}
(-b_j^i)\right). \]

\noindent
However, this condition is not sufficient to assure the Buchsbaum
property for $G$ as the following examples show.

\begin{example}
Consider the numerical semigroup of Example \ref{5,6,13} and
its Apery table:

\begin{table}[ht]
\renewcommand\arraystretch{1.5}
\[
\begin{array}{|l|c|c|c|c|c|}\hline \Ap(S)& 0&6&12&13&19\\ \hline
\Ap(M)&5&6&12&13&19\\ \hline \Ap(2M)&10&11&12&18&19\\ \hline
 \Ap(3M)&15&16&17&18&24\\ \hline
 \Ap(4M)&20&21&22&23&24 \\ \hline
 \end{array}.
 \]
 \end{table}
Then, we have that $0\neq (t^{6})^{\ast}, (t^{13})^{\ast} \in
\frac{\fm}{\fm^2} \subseteq G$, $(t^{13})^{\ast} \in T(G)$, and $0\neq (t^{6})^{\ast}\cdot
(t^{13})^{\ast}=\overline{t^{19}} \in \frac{\fm^2}{\fm^3}$ and so
$G$ is not Buchsbaum.
\end{example}

\begin{example}
\label{9,10,11,23} Let $S=\langle 9,10,11,23\rangle$. The Apery table is

\begin{table}[ht]
\renewcommand\arraystretch{1.5}
\[
\begin{array}{|l|c|c|c|c|c|c|c|c|c|}\hline \Ap(S)& 0&10&11&21&22&23&33&34&44 \\ \hline
\Ap(M)& 9&10&11&21&22&23&33&34&44 \\ \hline
 \Ap(2M)&18&19&20&21&22&32&33&34&44 \\ \hline
 \Ap(3M)&27&28&29&30&31&32&33&43&44 \\ \hline
 \Ap(4M) &36&37&38&39&40&41&42&43&44 \\ \hline
 \end{array}.
 \]
 \end{table}

Then, $0\neq (t^{11})^{\ast}, (t^{23})^{\ast} \in \frac{\fm}{\fm^2} \subseteq G$, $(t^{23})^{\ast} \in
T(G)$ and $0\neq (t^{11})^{\ast}\cdot (t^{23})^{\ast} \in \frac{\fm^2}{\fm^3}$ and so
$G$ is not Buchsbaum.
\end{example}

\begin{lemma}
\label{one box of length one} If $ \displaystyle{ G\cong F\oplus
\bigoplus_{i=1}^{e-1} F(-d_i)
 \oplus \frac{F}{(x^{\ast})F}
(-b)}$, then $G$ is Buchsbaum.
\end{lemma}
\begin{proof}
The statement is clear since, in this case, the torsion submodule $T(G)$ coincides with the socle of $G$.
\end{proof}

\begin{corollary}
Let $A$ be a numerical semigropup ring of multiplicity $4$. Then, its tangent is always Buchsbaum.
\end{corollary}

In other terms, the above lemma says that the tangent cone $G$ of
a numerical semigroup ring that verifies $\la (H^0_{G_{+}}(G))
\leq 1$ is Buchsbaum. Victoria A. Sapko has conjectured in
\cite{S} that the converse is true for the case of a 3-generated
semigroup ring. Recently, Yi Huang Shen \cite{Sh} has given a
positive answer to this conjecture, but on the basis of some of
our computations one could ask for a similar question for any
numerical semigroup ring.

\bibliographystyle{amsalpha}

\begin{thebibliography}{A}


\bibitem[BF]{BF} V. Barucci, R. Fr\"{o}berg,
\textit{Associated graded rings of one-dimensional analytically
irreducible rings}, J. Algebra \textbf{304} (2006), 349--358.

\bibitem[CZ1]{CZ1} T. Cortadellas, S. Zarzuela,
\textit{On the structure of the fiber cone of ideals with analytic
spread one},  J. Algebra  \textbf{317}  (2007),{no}. 2, 759--785.

\bibitem[CZ2]{CZ2} T. Cortadellas, S. Zarzuela,
\textit{Apery and microinvariants of a one dimensional Cohen-Macaulay local ring and the invariants of its tangent cone},  preprint, 2008.

\bibitem[E]{E} \textit{J. Elias, On the deep structure of the blowing-up of
curve singularities}, Math. Proc. Camb. Phil. Soc. \textbf{131}
(2001), 227--240.

\bibitem[GAP4]{GAP4}
The GAP Group,\textit{ GAP - Groups, Algorithms, and Programming-
Version 4.4.10} (2007) (http://www.gap-system.org).

\bibitem[NumericalSgps]{NumericalSgps}
M. Delgado, P. A. Garcia-S\'{a}nchez, J. Morais  ,
\textit{NumericalSgps - a GAP package, 0.95} (2006),
(http://www.gap-system.org/Packages/numericalsgps).

\bibitem[S]{S} V. A. Sapko, \textit{Associated graded rings of numerical semigroup
rings},  Comm. Algebra  29  (2001),  no. 10, 4759--4773.

\bibitem[Sh]{Sh} Y. Shen, \textit{Tangent cones of numerical semigroup rings with small
embedding dimension}, arXiv:0808.2162v1[math.AC].

\end{thebibliography}

\end{document}